\documentclass[11pt]{amsart}
 
\usepackage{amsfonts}
\usepackage{amssymb}
\usepackage{amscd}
\input xy
\xyoption{all}
 
\newcommand{\marginlabel}[1]%
  {\mbox{}\marginpar{\raggedleft\hspace{0pt}\bfseries\sf#1}}

\def\CC{{\mathbb C}}
\def\PP{{\textbf P}}
\def\OO{{\mathcal O}}

\def\F{\mathcal{F}}
\def\E{\mathcal{E}}
\def\G{\mathcal{G}}

\def\I{\mathcal{I}} 

\def\cM{\mathcal{M}}

\def\Pic0{{\rm Pic}^0(X)}

\theoremstyle{plain}

\newtheorem{theorem}{Theorem}[section]
\newtheorem{proposition}[theorem]{Proposition}
\newtheorem{corollary}[theorem]{Corollary}
\newtheorem{lemma}[theorem]{Lemma}

\theoremstyle{definition}
\newtheorem{definition}[theorem]{Definition}
\newtheorem{remark}[theorem]{Remark}
\newtheorem{example}[theorem]{Example}
   
\newtheorem{conjecture}[theorem]{Conjecture}
\newtheorem{conjecture/question}[theorem]{Conjecture/Question}
\newtheorem{question}[theorem]{Question}  
\newtheorem{remark/definition}[theorem]{Remark/Definition} 
 
\theoremstyle{remark}

\begin{document}
 
\title{M-regularity and the Fourier-Mukai transform}

\author[G. Pareschi]{Giuseppe Pareschi} 
\address{Dipartamento di Matematica, Universit\`a di Roma, Tor Vergata, V.le della 
Ricerca Scientifica, I-00133 Roma, Italy}
\email{{\tt pareschi@mat.uniroma2.it}}

\author[M. Popa]{Mihnea Popa}
\address{Department of Mathematics, University of Chicago,
5734 S. University Ave., Chicago, IL 60637, USA }
\email{{\tt mpopa@math.uchicago.edu}}

\thanks{2000\,\emph{Mathematics Subject Classification}.
 Primary 14K05; Secondary 14K12, 14E25, 14F05.
\newline MP was partially supported by the NSF grant DMS 0500985
and by an AMS Centennial Fellowship}

\maketitle

$${\rm Abstract}$$
This is a survey of $M$-regularity and its applications, expanding on lectures given by the second author at the Seattle conference, in August 2005, and at the Luminy workshop "G\' eom\' etrie Alg\' ebrique Complexe", in October 2005\footnote{He would like to thank A. Beauville for the invitation, and the organizers of the GAC workshop for providing this opportunity.}.

\tableofcontents

\markboth{G. PARESCHI and M. POPA}
{\bf M-REGULARITY AND THE FOURIER-MUKAI TRANSFORM}

\section{Global generation of coherent sheaves} 

In these notes we will mostly work with smooth projective varieties over an arbitrary 
algebraically closed field, although in many applications we will restrict to the case of
characteristic zero. 

Let $X$ be a smooth projective variety.  Many basic problems related to the 
geometry of $X$ can be expressed in terms of the global generation of appropriately 
chosen coherent sheaves on X. For instance, we usually ask whether a given linear 
series on $X$ gives a morphism, or better an embedding, and if that is the case, 
what are the equations cutting out the image of $X$ in projective space. 

\begin{example}\label{basic}
(1) Given a line bundle $A$, the question of whether the linear series $|A|$ is basepoint-free 
is tautologically a global generation question. Most commonly, given an ample line bundle 
$L$ on $X$, we wonder for what values of $m$ the linear series $|mL|$, $|K_X + mL|$, or 
even $|m(K_X + L)|$ are basepoint-free\footnote{For the first two, the existence of such $m$ 
follows by definition. For the last, we have to assume some sort of positivity for $K_X +L$ -- if it is nef, 
for instance, the existence of such an $m$ follows from Kawamata's Basepoint-Free theorem.}.  

(2) A line bundle $L$ on $X$ is very ample if and only if 
$L\otimes \I_{x}$ is globally generated for all $x\in X$. Thus the very ampleness problem can 
be reduced to global generation. 

(3) If $X\subset \PP^n$, then the global generation of the twisted ideal sheaf $\I_X (d)$ implies that 
$X$ is locally cut out by equations of degree $d$. Similarly, if $X$ is a subvariety of a polarized 
variety $(Y, L)$, the global generation of the sheaf $\I_X \otimes L^d$ gives information about the 
equations describing locally $X$ inside $Y$. 
\end{example}

The question then becomes how to check global generation of coherent sheaves effectively. For sheaves 
on $\PP^n$, an important technique in this direction was introduced by Mumford \cite{mumford3}, and   
later called \emph{Castelnuovo-Mumford regularity}. A coherent sheaf $\F$ on $\PP^n$ is called $m$-regular in the sense of Castelnuovo-Mumford if 
$$H^i (\F (m-i)) = 0, {\rm ~for~all~} i>0.$$ 
The main result involving this concept is the following: 

\begin{theorem}[Castelnuovo-Mumford Lemma]\label{cm_lemma}
Let $\F$ be a $0$-regular coherent sheaf on $\PP^n$. Then the following hold: 
\begin{enumerate}
\item $\F$ is globally generated. 
\item $\F$ is $m$-regular for all $m \ge 1$.
\item The multiplication map 
$$H^0 (\F) \otimes H^0 (\OO_{\PP^n}(k)) \longrightarrow H^0 (\F (k))$$ 
is surjective for all $k\ge 0$.
\end{enumerate}
\end{theorem}

As it is often possible to check the vanishing of cohomology, this result is one of the most widely used 
tools for attacking the problems mentioned above -- see \cite{positivity} I.1.8 for a detailed discussion 
and applications.  

\bigskip 

One of the key observations of \cite{pareschi}, \cite{pp1} and \cite{pp2}  is that on abelian varieties essentially all geometric problems  related to linear series and defining equations can be reduced to the global generation of suitable coherent sheaves. This goes via the notion of skew-Pontrjagin product, an operation similar to basic operations performed on cohomology classes on complex tori, and via surjectivity for multiplication maps on global sections. 

\begin{definition}
Let $X$ be an abelian variety of dimension $g$. 
Given two coherent sheaves $\E$ and $\G$ on $X$, 
their \emph{skew-Pontrjagin product} (see \cite{pareschi} \S1) is defined as 
$$\E\hat * \G:= d_*(p_1^*\E \otimes p_2^*\G), $$
where $d: X \times X \rightarrow X$ is the difference map $(x,y)\to x-y$. 
\end{definition} 

We will also use the following notation: given two sheaves $\E$ and $\G$ on $X$, 
we denote by $\cM( \E, \G)$ the locus of $x\in X$ where the multiplication map 
$$m_x:H^0(t_x^*\E)\otimes H^0(\G)\rightarrow H^0((t_x^*\E)\otimes \G)$$
is not surjective. (Here $t_x : X \rightarrow X$ denotes the morphism given by 
translation by $x$.) The relationship between skew Pontrjagin products and 
multiplication maps is provided by the following:

\begin{proposition}[\cite{pareschi} Proposition 1.1]\label{pontrjagin}
Let $\E$ and $\G$ be sheaves on $X$ such that $$H^i (t_x^* \E \otimes \G) = 0, 
{\rm ~for~ all~}i > 0.$$ Then $\cM(\E,\G)$ is precisely the locus where $ \E \hat * \G $ 
is not globally generated.   
\end{proposition}

In other words, if the skew-Pontrjagin product $\E \hat * \G $ is globally generated, then 
all the multiplication maps $m_x$ are surjective. In particular this is the case for 
$$ H^0(\E) \otimes H^0(\G) \longrightarrow H^0( \E \otimes \G).$$ 
It is to such multiplication maps that the geometric problems mentioned above are reduced. 

\begin{example}
(1) Let's assume that the line bundle $L$ on $X$ is very ample. In this case, the \emph{projective
normality} of $X$ in the embedding given by $L$ is known to be equivalent to the surjectivity of 
$$H^0(L) \otimes H^0 (L^k) \longrightarrow H^0( L^{k+1})$$ 
for all $k\ge1$. As ample line bundles on abelian varieties have no higher cohomology, by the 
remarks above this is implied by the global generation of the skew-Pontrjagin products 
$L \hat * L^k$. 

(2) Going one step further, consider the vector bundle $M_L$ defined as the kernel of the evaluation 
map on global sections of $L$: 
$$0\rightarrow M_L\rightarrow H^0(L)\otimes {\mathcal O}_X\rightarrow
L\rightarrow 0.$$ 
The global generation of the skew-Pontrjagin product 
$L\hat * (M_L\otimes L)$ ensures that  the homogeneous ideal of $X$ in the embedding given by $L$ is generated 
by quadrics. We will see this in \S4, together with a similar -- but gradually more complicated -- approach that can be adopted towards understanding the syzygies of $X$ in this embedding.  
\end{example} 

It becomes thus important to have criteria guaranteeing the global generation of coherent sheaves 
on abelian varieties. Theorem \ref{cm_lemma} serves as a good guideline, but we will see that on 
abelian varieties there is a more efficient notion of regularity. We will describe this in \S3, after a brief 
reminder of the Fourier-Mukai transform in \S2.

\section{The Fourier-Mukai transform}

In what follows, unless otherwise specified, $X$ will be an abelian variety 
of dimension $g$ over an algebraically closed field. We denote by $\hat{X}$ the dual abelian variety, 
which will often be identified with $\Pic0$. By $\mathcal{P}$ we denote a Poincar\'e
line bundle on $X\times \hat{X}$, normalized such that 
$\mathcal{P}_{|X\times\{0\}}$ and $\mathcal{P}_{| \{0\}\times \widehat{X}}$ are trivial.

We briefly recall the Fourier-Mukai setting, referring to \cite{mukai} for details.
To any coherent sheaf $\F$ 
on $X$ we can associate the sheaf ${p_{2}}_{*}({p_{1}}^{*}\F\otimes\mathcal{P})$ on
$\widehat{X}$, where $p_1$ and $p_2$ are the natural projections on $X$ and $\hat{X}$. 
This correspondence gives a functor
$$\hat{\mathcal{S}}: {\rm Coh(X)}\rightarrow {\rm Coh}(\widehat{X}).$$ 
If we denote by ${\bf D}(X)$ and 
${\bf D}(\hat{X})$ the (bounded) derived categories of Coh$(X)$ and Coh$(\hat{X})$, then 
the derived functor $\mathbf{R}\hat{\mathcal{S}}:{\bf D}(X)\rightarrow {\bf D}(\hat{X})$ 
is defined and called the \emph{Fourier-Mukai functor},
and one can consider $\mathbf{R}\mathcal{S}:{\bf D}(\hat{X})\rightarrow {\bf D}(X)$ 
in a similar way. Mukai's main result is the following: 

\begin{theorem}[\cite{mukai}, Theorem 2.2]\label{duality}
The Fourier-Mukai functor establishes an equivalence of categories between ${\bf D}(X)$ and 
${\bf D}(\widehat{X})$. More precisely there are isomorphisms of functors:
$$\mathbf{R}\mathcal{S}\circ\mathbf{R}\widehat{\mathcal{S}}\cong (-1_X)^{*}
[-g] ~{\rm ~and~}~ \mathbf{R}\widehat{\mathcal{S}}\circ\mathbf{R}\mathcal{S}\cong (-1_{\widehat{X}})^{*} [-g].$$
\end{theorem}

Most important to us is the cohomological information encoded by this theorem:  

\begin{corollary}\label{ext}
If $A$ and $B$ are objects in ${\bf D}(X)$, then 
$${\rm Ext}^i_{{\bf D}(X)} (A, B) \cong {\rm Ext}^i_{{\bf D}(\widehat{X})} (\mathbf{R}\hat{\mathcal{S}}A, 
\mathbf{R}\hat{\mathcal{S}}B).$$
\end{corollary}

One more piece of terminology is useful in what follows. 
\begin{definition}
Let $R^j\hat{\mathcal{S}}(\F)$ be the cohomologies of the derived complex 
$\mathbf{R}\hat{\mathcal{S}}(\F)$. Following \cite{mukai}, 
the sheaf $\F$ satisfies W.I.T. (the \emph{weak index theorem}) with index $i$ if 
$$R^j\hat{\mathcal{S}}(\F)=0, {\rm ~for~ all~} j\neq i .$$ 
It satisfies the stronger I.T. (the
\emph{index theorem}) with index $i$ if 
$$H^{i}(\F\otimes\alpha)=0, {\rm~ for~ all~} \alpha\in {\rm Pic}^{0}(X) {\rm~ and~ all~} i\neq j.$$
By the Base Change theorem, in this last situation $R^{j}\hat{\mathcal{S}}(\F)$ is locally free. 
If $\F$ satisfies W.I.T. with index $i$, $R^{i}\hat{\mathcal{S}}(\F)$ is denoted by 
$\hat{\F}$ and called the $\it{Fourier~
transform}$ of $\F$. Note that then $\mathbf{R}\hat{\mathcal{S}}(\F)\cong \hat{\F}[-i]$.
\end{definition}

\section{The main results on M-regularity} 

Let $X$ be an abelian variety of dimension $g$.

\begin{definition}
A coherent sheaf $\F$ on $X$ is called \emph{$M$-regular} (or  
\emph{$Mukai$-regular}) if
$${\rm codim}_{\hat X}({\rm Supp}(R^i\hat{\mathcal S}(\F)))>i {\rm ~ for ~all~} i=1,\ldots ,g$$ 
(where, for $i=g$, this means that ${\rm Supp}(R^g\hat{\mathcal S}(\F)) = \emptyset$).
\end{definition}

\begin{remark}
By Base Change we see that there is always an inclusion
${\rm Supp}(R^i\hat{\mathcal S}(\F))\subset V^i(\F)$, where $V^i(\F)$ is the \emph{cohomological 
support locus} (cf. 
\cite{gl}):
$$V^i(\F):= \{ \alpha ~|~h^i(\F\otimes \alpha )\neq 0\}\subset {\rm Pic}^0(X).$$
Consequently, $M$-regularity is achieved if in particular
$${\rm codim}(V^i(\F))>i {\rm ~for~ all~} i=1,\ldots ,g.$$
It is this property that one usually checks in applications.
\end{remark}

The main theme of these notes is that $M$-regularity is precisely the cohomological 
condition to be checked in order to obtain geometric results via global generation. Here 
is the main result: 

\begin{theorem}[\cite{pp1}, Theorem 2.4]\label{mreg}
Let $\F$ be a coherent sheaf and $L$ an invertible sheaf supported on a 
subvariety $Y$ of the abelian variety $X$ (possibly $X$ itself). If both $\F$ and $L$ are 
$M$-regular as sheaves on $X$, then $\F\otimes L$ is globally generated.
\end{theorem}

\begin{example}
It is worthwhile noting that the statement does not hold if the hypothesis that $L$ be a line 
bundle is dropped, even if it is locally free. Consider for instance a principally polarized abelian 
variety $(X, \Theta)$. (For simplicity we identify it with its principally polarized dual via the polarization). Define 
$$F: = \OO_X(\Theta) \otimes \widehat{\OO_X(n\Theta)} {\rm ~and ~ } L:= \widehat{\OO_X(-n\Theta)}.$$
for some integer $n > 1$. By \cite{mukai} 3.11 we have that $\phi_n^* \widehat{\OO_X(-n\Theta)} \cong \bigoplus \OO(n\Theta)$, where $\phi_n$ denotes multiplication by $n$ on $X$. This implies that both 
$F$ and $L$ satisfy I.T. with index $0$, the strongest form of $M$-regularity. However, by applying the 
trace map to $\widehat{\OO_X(n\Theta)} \otimes \widehat{\OO_X(-n\Theta)}$ we see that $\OO_X(\Theta)$ is a direct summand of $F\otimes L$, so the latter cannot be globally generated. 
\end{example}

We will describe briefly the ideas involved in the proof of Theorem \ref{mreg}. Some of the 
intermediate results will also have independent applications. The main point is to consider a modified notion of generation by global sections. 

\begin{definition}\label{cgg}
Let $Y$ be a  variety. A coherent sheaf $\F$ on $Y$ is  
\emph{continuously globally generated} if for any non-empty open subset 
$U\subset {\rm Pic}^{0}(Y)$ the sum of evaluation maps 
$$\bigoplus_{\alpha\in U} H^0(\F\otimes \alpha)\otimes \alpha^\vee \longrightarrow \F$$ 
is surjective. 
\end{definition}

Theorem \ref{mreg} follows immediately from the fact that the tensor product of a continuously globally generated sheaf and a continuously globally generated line bundle is globally generated\footnote{This is seen by simply mimicking the use of translations to see that the tensor product of two ample line bundles has no base locus.}, once 
we have the following result\footnote{Though not presented as such because of the new terminology, this should in fact be considered the main result of the section.}. 

\begin{proposition}[\cite{pp1}, Proposition 2.13]\label{cont}
Any $M$-regular coherent sheaf is continuously globally generated. 
\end{proposition}

This result follows in turn in a standard way from a technical statement which generalizes ideas that have already appeared in a more restrictive setting in  work of Mumford, Kempf and Lazarsfeld. 

\begin{theorem}[\cite{pp1}, Theorem 2.5]\label{amultiplication}
Let $\F$ and $H$ be sheaves on $X$ such that $\F$ is $M$-regular and 
$H$ is locally free satisfying IT with index $0$. Then, for any non-empty Zariski open set 
$U\subset \Pic0$, the map
$$\bigoplus_{\alpha\in U}H^0(X,\F\otimes \alpha)\otimes H^0(X,H\otimes \alpha^{-1})
\buildrel{\oplus m_{\alpha}}\over
\longrightarrow H^0(X,\F\otimes H)$$
is surjective, where $m_{\xi}$ denote the multiplication maps on global sections.
\end{theorem}

It is here that Mukai's results on the Fourier functor are heavily used. Let's give a glimpse of how it comes into the picture. For example, by Serre duality followed by Mukai's duality, more precisely 
Corollary \ref{ext}, we have that 
$$H^0 (\F \otimes H)^{\vee} \cong 
{\rm Ext}^g_{{\bf D}(X)} (\F, H^{\vee}) \cong {\rm Ext}^g_{{\bf D}(\widehat{X})} (\mathbf{R}\hat{\mathcal{S}} \F,  \mathbf{R}\hat{\mathcal{S}} H^{\vee}).$$
Note that $H^\vee$ satisfies W.I.T. with index $g$, and so $\mathbf{R}\hat{\mathcal{S}} H^{\vee} \cong 
\widehat{H^\vee}[-g]$. Via a geometric reduction, the statement of the theorem eventually follows from a study of the standard spectral sequence 
$$E^{ij}_2 ={\rm Ext}^i(R^j\hat{\mathcal{S}} \F,\widehat{H^\vee})\Rightarrow
{\rm Ext}^{i-j}_{{\bf D}(\hat X)}({\bf R}\hat{\mathcal{S}} \F,\widehat{H^\vee}),$$
where the $M$-regularity of $\F$ ensures by definition the vanishing of $E^{ij}_2$ for $0 < i \le j$. 
It is worth noting that there is a "preservation of vanishing" statement, important for deducing the index theorem from $M$-regularity.  

\begin{proposition}[\cite{pp1}, Proposition 2.9]\label{freg_vanishing}
Let $\F$ be an $M$-regular coherent sheaf on $X$ and $H$ a locally free sheaf satisfying I.T. with 
index $0$. Then $\F\otimes H$ satisfies I.T. with index $0$.
\end{proposition} 

The following natural question is still open. We expect a positive answer, at least when one of the 
two terms is locally free. 

\begin{question}
\emph{Is the tensor product of two $M$-regular sheaves $M$-regular?} 
\end{question}

\noindent
\textbf{First examples.}
To have a first look at how this works, here are some very basic examples. 

\begin{example}\label{first}
(1) For a line bundle $L$ on an abelian variety $X$:
$$L {\rm~is ~M-regular} \iff L {\rm ~satisfies~ I.T. ~with~index~}0 \iff L {\rm ~is~ample}.$$
Indeed, ample line bundles on abelian varieties have no higher cohomology, so the 
two implications from right to left are immediate. On the other hand, if $L$ is $M$-regular, by definition and Proposition \ref{cont} we see that for general $\alpha\in \Pic0$ we have 
$h^i (L\otimes \alpha) = 0$ for $i>0$, and $h^0 (L\otimes \alpha) \neq 0$. In Mumford's terminology (cf. \cite{mumford} \S16) this means that $L\otimes \alpha$ is a non-degenerate line bundle of index $0$, which is equivalent to its ampleness  (cf. \emph{loc. cit.}, p.60).
This in turn is equivalent to the ampleness of $L$ itself. (Cf. Lemma \ref{nonvanishing} and 
Proposition \ref{ampleness} for much more general statements.)

(2) A line bundle $L$ on a smooth curve $C$ of genus $g\geq 1$ is $M$-regular (via an 
Abel-Jacobi embedding $C\subset J(C)$) if and only if $d={\rm deg}~L \geq g$.
Indeed, the $M$-regularity of $L$ is equivalent by Riemann-Roch and base-change to the fact that
the Brill-Noether locus $W_d^{d-g+1}$ has codimension at least $2$ in ${\rm Pic}^d(C)$. But this is easily seen to be equivalent to $d\geq g$. A similar result holds for the image of a box-product of $d$ 
line bundles via the desymmetrized Abel-Jacobi mapping $C^d \rightarrow W_d \subset J(C)$ (cf. \cite{pp1} Example 3.2). 

(3) If $L$ is an ample line bundle on the abelian variety $X$, then: 
$$L\otimes \I_{\{x\}} {\rm~is ~M-regular},~ \forall~x\in X \iff {\rm codim}~ {\rm Bs}|L| \ge 2, $$ 
i.e. if $|L|$ has no base divisor.  This follows from the standard cohomology sequence 
$$0 \longrightarrow L\otimes \I_{\{x\}} \longrightarrow L \rightarrow L\otimes \OO_x \longrightarrow 
0$$ by passing to cohomology, combined with the fact that twisting $L$ with line bundles in ${\rm Pic}^0(X)$ is equivalent to translating.
\end{example} 

Applied to these examples, Theorem \ref{mreg}, together with Example \ref{basic} (2), implies the 
following well-known starting points in the study of linear series on curves and abelian varieties. 

\begin{corollary}
Let $X$ be an abelian variety  and $C$ a smooth projective curve. Then:
\newline
\noindent
(i) Let $L$ be a line bundle of degree $d$ on $C$. If $d \ge 2g$, then 
$L$ is globally generated, and if $d\ge 2g+1$, then $L$ is very ample.
\newline
\noindent
(ii)(\emph{Lefschetz Theorem}) If $L$ is an ample line bundle
on $X$, then $L^2$ is globally generated and  $L^3$ is very ample.
\newline
\noindent
(iii)(\emph{Ohbuchi's Theorem}) If $L$ is an ample line bundle on $X$ 
with no base divisor, then $L^2$ is very ample.
\end{corollary}

\noindent
\textbf{Theta regularity.} 
$M$-regularity is not an immediately obvious analogue of the usual
notion of Castelnuovo-Mumford regularity on $\PP^n$. 
There is however a particular instance of the general theory, depending on 
a fixed polarization $\Theta$, where one clearly sees a similarity -- this is called 
\emph{Theta regularity} in \cite{pp1} and \cite{pp3}.

\begin{definition}
A coherent sheaf $\F$ on a polarized abelian variety $(X, \Theta)$ is 
called $m$-$\Theta$-\emph{regular} 
if $\F((m-1)\Theta)$ is $M$-regular.
If $\F$ is $0$-$\Theta$-regular, we will simply call it $\Theta$-\emph{regular}.
\end{definition}

\begin{example}
If $C$ is a smooth projective curve and $W_d \subset J(C)$ is the image of the 
$d$-th symmetric product of $C$ in the Jacobian, then 
 $\OO_{W_d}$ is $2$-$\Theta$-regular and $\I_{W_d}$ is $3$-$\Theta$-regular
 (cf. \S6, Theorem \ref{w_d} and its proof). 
 \end{example}

With this language we have an abelian analogue of the Castelnuovo-Mumford Lemma:

\begin{theorem}[\cite{pp1}, Theorem 6.3]\label{acm}
Let $\F$ be a $\Theta$-regular coherent sheaf on $X$. Then:
\newline
\noindent
(1) $\F$ is globally generated.
\newline
\noindent
(2) $\F$ is $m$-$\Theta$-regular for any $m\geq 1$.
\newline
\noindent
(3) The multiplication map
$$H^0(\F(\Theta))\otimes H^0(\OO(k\Theta))\longrightarrow H^0(\F((k+1)\Theta))$$
is surjective for any $k\geq 2$.
\end{theorem}

\begin{remark}
The statement in (3) is optimal, as it follows 
for example by considering $\F$ equal to $\OO(2\Theta)$ (when we cannot make 
$k=1$). This particular example is the 
well-known projective normality result for multiples of ample line 
bundles treated by Koizumi \cite{koizumi} and Mumford  \cite{mumford2}.
\end{remark}

\noindent
{\bf WIT-regularity.} This is a very useful alternative to $M$-regularity -- or, maybe 
more correctly said, to the Index Theorem with index $0$ condition -- when we are not 
looking for full global generation, but rather just a good control of the locus where the 
sheaf in question is not globally generated. Its use is most significant in the syzygy questions described in \S4. 

We can weaken the condition of continuous global generation (Definition \ref{cgg}) as follows (cf. \cite{pp2}, Definition 2.1): 
 
\begin{definition}
Given a sheaf $\F$, we define its \emph{Fourier jump locus} as the 
locus $J(\F)\subset \Pic0$ consisting of $\alpha\in \Pic0$ where
$h^0(\F\otimes \alpha)$ jumps, i.e. it is different from its minimal value over $\Pic0$. 
$\F$ is said to be \emph{weakly continuously generated} if the map
$$\bigoplus_{\alpha\in U}H^0(\F\otimes \alpha)\otimes
\alpha^{-1} \longrightarrow \F$$ 
is surjective for \emph{any non-empty Zariski-open set $U\subset \hat X$ containing
$J(\F)$}. 
\end{definition}

\begin{theorem}[\cite{pp2}, Theorem 4.1]\label{wit}
Let $F$ be a locally free sheaf on $X$ such that
$F^\vee$ satisfies W.I.T. with index $g$ and the torsion part of
$\widehat{F^\vee}$ is a torsion-free sheaf on a reduced subscheme $Y$ of
$X$. Then the following hold:
\newline
\noindent
(a) $F$ is weakly continuously generated.
\newline
\noindent
(b) Let morever $A$ be a continuously globally generated line bundle on $Y$. Then 
\begin{itemize}
\item[(i)] $F\otimes A$ is generically globally generated.
\item[(ii)] If the Fourier-jump locus $J(F)$ is finite then 
$$B(F \otimes A)\subset \bigcup_{\xi\in J(F)}B(A\otimes P_\xi),$$ 
where $B(\cdot)$ denotes the locus where the sheaf is not globally generated. 
\end{itemize}
\end{theorem}

The underlying point here is that locally free sheaves satisfying I.T. with index $0$ have duals satisfying the W.I.T. with index $g$, but the converse is not true, and so Theorem
\ref{mreg} does not necessarily apply. 

\begin{example} 
Typical such behavior is exemplified already by the structure sheaf $\OO_X$, which satisfies W.I.T. with index $g$. It is not $M$-regular, and its dual also satisfies W.I.T. with index $g$, since it is the same sheaf. The Fourier jump locus is the point-set $\{0\}$ on $\hat X$, which obviously has to appear among the elements of $U$ in order to have continuous generation. 
\end{example} 

\noindent
\textbf{Generic Vanishing.} 
If $F$ is a locally free sheaf on $X$ satisfying the hypothesis of Theorem \ref{wit}, it is explained by 
Hacon in \cite{hacon} that $F$ satisfies a Green-Lazarsfeld-type \cite{gl} Generic Vanishing condition, namely 
$${\rm codim}( V^i (F)) \ge i, {\rm~for~all~} i.$$ 
(One can show that this is in fact equivalent to the condition that ${\rm codim}~{\rm Supp}
(R^i \hat{\mathcal{S}} F) \ge i$ for all $i$ -- thus it is uniformly one step ``worse" than $M$-regularity.) 

In  the forthcoming \cite{pp5} it is shown that 
the converse is also true. Thus a posteriori weak continuous generation is essentially a property naturally associated to sheaves satisfying Generic Vanishing, just like continuous generation is associated to the stronger $M$-regularity.  We also show that, if a locally free sheaf $F$ satisfies  Generic Vanishing, then $F$ is $M$-regular if and only if $\widehat{F^\vee}$ is \emph{torsion-free}.
One obtains a natural commutative algebra interpretation of the $M$-regularity condition. 
In brief, we have that 
$$R^i \hat{\mathcal{S}} F \cong \mathcal{E} xt^i (\widehat{F^\vee}, \OO_X) 
{\rm ~ for~all~} i,$$ 
and the support of the $i$-th such ${\rm Ext}$-sheaf of a torsion-free sheaf always has codimension greater 
than $i$ by the Syzygy Theorem.  
These results are particular consequences -- in the context of abelian varieties -- of 
a systematic analysis of Generic Vanishing conditions with respect to arbitrary Fourier-Mukai functors 
on projective varieties. This is the main theme of \cite{pp5}.

\section{Equations and syzygies of abelian varieties} 

To every ample line bundle on an abelian variety one can associate an invariant 
called the $M$-regularity index, defined below. This invariant governs 
the higher order properties of embeddings of abelian varieties by powers of ample 
line bundles, and also the equations and the syzygies of such embeddings (here though for the most part conjecturally). For high values, the $M$-regularity index is still quite 
mysterious, however for small values it has nice geometric interpretation. 
 
\begin{definition}[\cite{pp2}, \S3]
The $M$-\emph{regularity index} of an ample line bundle $A$ is defined as 
$$m(A):={\rm max}\{l~|~A\otimes 
m_{x_1}^{k_1}\otimes \ldots \otimes m_{x_p}^{k_p}~{\rm is ~}M{\rm -regular~ for 
~all~distinct~}$$
$$x_1,\ldots,x_p\in X {\rm~with~} \Sigma k_i=l\}.$$
\end{definition}

\begin{example}{\bf (Small values of $m(A)$.)}\label{base_div}
By Example \ref{first}(1), having $m(A)\ge 0$ is equivalent to the ampleness of $A$, with no restrictions. Furthermore, by Example \ref{first}(3), $m(A)\geq 1$ if and only if $A$ 
does not have a base divisor. One can similarly see that if $A$ gives a birational map which 
is an isomorphism outside a codimension $2$ subset, then $m(A)\geq 2$. 
\end{example}

\begin{definition}
A line bundle $L$ is called $k$-\emph{jet ample}, $k\geq 0$, if the restriction map
$$H^0(L)\longrightarrow H^0(L\otimes \OO_X/
m_{x_1}^{k_1}\otimes \ldots \otimes m_{x_p}^{k_p})$$ 
is surjective for any distinct points $x_1,\ldots,x_p$ on $X$ such that $\Sigma k_i =k+1$.
(In particular $0$-jet ample means 
globally generated, $1$-jet ample means very ample.)
\end{definition}

A natural generalization of the classical theorem of Lefschetz and of its extensions by 
Ohbuchi \cite{ohbuchi} and Bauer-Szemberg \cite{bs} is the following: 

\begin{theorem}[\cite{pp2}, Theorem 3.8]
If $A$ and $M_1,\ldots ,M_{k+1-m(A)}$ are ample line bundles on $X$, $k\geq m(A)$, 
then $A\otimes M_1\otimes \ldots\otimes M_{k+1-m(A)}$ is $k$-jet ample. In particular
$A^{\otimes(k+2-m(A))}$ is $k$-jet ample.
\end{theorem}

A simple corollary of this result is a lower bound for the Seshadri constant of an ample 
line bundle on an abelian variety in terms of the $M$-regularity index and its 
asymptotic version. 

\begin{corollary}[\cite{pp3}, Theorem 3.4]
If $L$ is an ample line bundle on the abelian variety $X$ and $\epsilon(L)$ denotes its Seshadri constant\footnote{The definition and properties of Seshadri constants are explained for example 
in \cite{positivity}, Ch.V. The particular case of abelian varieties is treated in \emph{loc. cit.} \S5.3.}, 
then 
$$\epsilon(L) \ge \rho (L) : = \underset{n}{{\rm sup}} \frac{m(L^n)}{n} \ge {\rm max}
\{ m(L), 1\}.$$
\end{corollary}

\noindent
{\bf Equations and syzygies.}
Given a variety $X$ embedded in projective space by a complete linear
series $|L|$, the line bundle $L$ is said to \emph{satisfy property
$N_p$ } if the first $p$ steps of the minimal graded free resolution of the
algebra
$R_L=\bigoplus H^0(L^n)$ over the polynomial ring $S_L=\bigoplus {\rm Sym}^nH^0(L)$
are linear, i.e. of the form
$$  S_L(-p-1)^{\oplus i_{p}}\rightarrow 
S_L(-p)^{\oplus i_{p-1}}\rightarrow\cdots\rightarrow 
S_L(-2)^{\oplus i_1}\rightarrow S_L\rightarrow R_L\rightarrow 0.$$
Thus $N_0$ means that the embedded variety is projectively normal
(\emph{normal generation} in Mumford's terminology), $N_1$ means that
the homogeneous ideal is generated by quadrics (\emph{normal
presentation}), $N_2$ means that the relations among these quadrics are
generated by \emph{linear} ones and so on. 

Koizumi \cite{koizumi} proved that in the first embedding suggested by the Lefschetz
theorem, namely that given by $L^3$, the variety is projectively normal, i.e. it satisfies $N_0$. 
In Mumford's fundamental work on the equations defining abelian varieties \cite{mumford2} and on varieties cut out by quadratic equations \cite{mumford4}, and in improvements by 
Kempf \cite{kempf1}, it is shown that in the embedding given by $L^4$ the variety $X$ is cut out by quadrics, i.e. it satisfies $N_1$. The first author proved in \cite{pareschi} the following conjecture of Lazarsfeld, improving on further work by Kempf \cite{kempf2} and generalizing the above results to arbitrary syzygies.  

\begin{theorem}[\cite{pareschi}, Theorem 4.3] 
Let $X$ be an abelian variety and $A$ a line bundle on $X$. 
If ${\rm char}(k)$ does not divide $(p+1)$ and
$(p+2)$, then for $k\ge p+3$ the line bundle $A^{k}$ satisfies property $N_p$.
\end{theorem} 

As in the well-known Green-Lazarsfeld picture for curves, this is not the end of the 
story, and finer regularity properties of the line bundle should imply better and better 
syzygy properties. In \cite{pp2} the following conjecture is formulated: 

\begin{conjecture}
\emph{  Let $p\ge m$ be non-negative integers.
If $A$ is ample and $m(A)\ge m$, then $A^{\otimes k}$ satisfies
$N_p$ for any $k\ge p+3-m$.} 
\end{conjecture}

The main result of \cite{pp2} is a proof of the conjecture above in the case $m(A)=1$, 
i.e. the case when $|A|$ has no base divisor -- this is the generic behavior of an ample 
line bundle on an irreducible abelian variety, and it is well understood that in this situation 
one should already have better results than in the case of arbitrary ample line bundles. 

\begin{theorem}[\cite{pp2}, Theorem 6.2]\label{syzygies} 
In the previous setting, assume in addition that the linear 
system $|A|$ has no base divisor. Then for $k\ge p+2$ the 
line bundle $A^{k}$ satisfies property $N_p$.
\end{theorem}

We give a rough sketch of the proof of the theorem in the $N_1$ case -- this already contains all of the main ideas involved. The statement is simply the following: 
\begin{itemize}
\item \emph{If $L$ is an ample line bundle with no base divisor, then the ideal of $X$ in the embedding given by $L^3$ is generated by quadrics.}
\end{itemize}

Denote $A := L^3$. We have the usual evaluation sequence:
\begin{equation}\label{evaluation}
0\rightarrow M_A\rightarrow H^0(A)\otimes {\mathcal O}_X\rightarrow
A\rightarrow 0.
\end{equation}
Since $H^1(A) =0$, by work of Green and Lazarsfeld (cf. e.g. \cite{lazarsfeld}) it is well known that property $N_1$ is implied by the vanishing $H^1 (\wedge^2 M_A \otimes A^k) = 0$ for all $k \ge 1$. We will further 
restrict to the hardest case, namely $k=1$. In addition, since we are in characteristic different from $2$, the vanishing we're after is implied by the stronger 
$$H^1 (\otimes^2 M_A \otimes A) = 0.$$ 
By twisting the sequence above by $M_A\otimes A$ and passing to cohomology, this is 
in turn equivalent to the surjectivity of the multiplication map
$$H^0(A)\otimes H^0 (M_A\otimes A) \longrightarrow H^0(M_A\otimes A^2).$$
By the discussion in \S1, this surjectivity (and more) is implied by the global generation of the skew-Pontrjagin product $A\hat * (M_A\otimes A)$. 

This is where $M$-regularity comes into play. We use the hypothesis on the base locus of $|A|$ in order to prove the following:

\noindent
{\bf Claim.} The sheaf  $[A\hat * (M_A\otimes A)]\otimes L^{-1}$ is $M$-regular. 

Theorem \ref{mreg} implies then that the skew-Pontrjagin product is 
globally generated, and so our result. For the Claim, it is enough to prove that 
the cohomological support loci
$$V^i:=\{~\alpha \in \Pic0 \> | \> h^i((A\hat *( M_{A}\otimes A))\otimes
L^{-1} \otimes \alpha)>0\}$$
have {\rm codim}ension $> i$ for all $i>0$. A  simple result of "exchange of Pontrjagin and tensor product" via the base change theorem (\cite{pp2}, Proposition 5.5(b)) implies that the support loci above are the same as
$$V^i=\{\alpha \in \Pic0 \> |\> h^i((A\hat *(L^{-1} \otimes \alpha))\otimes
M_{A}\otimes A)>0\}.$$
After twisting the evaluation sequence (\ref{evaluation}) by $(A\hat *(L^{-1} \otimes \alpha))\otimes A$ and passing to cohomology, it is easily seen that $V^i = \emptyset$
for $i\ge2$, and one is left with showing that ${\rm codim}~V^1 \ge 2$. 
The same exact sequence (\ref{evaluation}) shows that 
\begin{itemize}
\item $V^1$ coincides with the locus where the
multiplication map
$$H^0(A)\otimes H^0(
(A\hat *(L^{-1} \otimes \alpha)) \otimes A)\rightarrow
H^0((A\hat *(L^{-1} \otimes \alpha)) \otimes A^2)$$
is not surjective.
\end{itemize}
Applying Proposition \ref{pontrjagin} one more time, it is enough to show that the 
locus 
$$B\bigl(A\hat * ((A\hat *(L^{-1}
\otimes \alpha))\otimes A)\bigr)$$ 
where the double skew-Pontrjagin product in the 
brackets is not globally generated has codimension at least $2$. For this one needs to 
use a version of the main global generation criterion which allows for base loci.  This 
is provided by Theorem \ref{wit}. A combination of this and explicit calculus with Pontrjagin 
products shows eventually (cf. \cite{pp2}, p.17--18) that in fact\footnote{up to some harmless 
translation which we skip for simplicity} 
$$B\bigl(A\hat * ((A\hat *(L^{-1}
\otimes \alpha))\otimes A)\bigr) \subset \bigcup_{\eta\in {\hat
X}_3}B(L\otimes \eta \otimes \alpha^{-1}),$$
where $ {\hat X}_3$ denotes the locus of $3$-torsion points on the dual of $X$. The initial hypothesis on the codimension of the base locus of $L$ implies then what is needed.

\section{Vanishing theorems and varieties of maximal Albanese dimension} 

In this section we present some applications to the study 
of linear series on varieties with generically finite Albanese maps. The main 
tool is the use of vanishing theorems as input for $M$-regularity. 
We start with a simple remark. 

\begin{lemma}\label{nonvanishing}
Let $\F$ be a non-zero coherent sheaf on an abelian variety $X$. 
If $\F$ is $M$-regular, then it is nef \footnote{Recall that to each coherent sheaf $\F$ one can associate the scheme $\PP(\F) : = {\rm Proj}(\oplus_m {\rm Sym}^m (\F))$ over $X$, together with an invertible sheaf $\OO_{\PP(\F)}(1)$. Then $\F$ is called nef if $\OO_{\PP(\F)}(1)$ is so.},  $h^0 (\F) > 0$, and $\chi (\F) > 0$.  
\end{lemma} 
\begin{proof} 
By Proposition \ref{cont}, $\F$ is continuously globally generated. Since any $\alpha\in \Pic0$ is 
nef, this means $\F$ is a quotient of a nef vector bundle, so it is itself nef. It also 
means that for $\alpha\in \Pic0$ general we have $H^0 (\F\otimes \alpha) \neq 0$, which implies 
$h^0 (\F) > 0$ by semicontinuity. 
By the definition of $M$-regularity, we have 
$$H^i (\F \otimes \alpha) = 0, {\rm ~for ~ all~} i>0 {\rm~and~} \alpha \in \Pic0 
{\rm ~general}.$$
Thus for $\alpha\in \Pic0$ general we have 
$$\chi (\F) = \chi (\F \otimes \alpha) =  h^0 (\F \otimes \alpha) > 0,$$ 
as the Euler characteristic is invariant under deformation. 
\end{proof}

\begin{remark}
The nefness of $\F$ in the Lemma is easily proved, but is rather weak. 
Debarre \cite{debarre2} has in fact shown that every $M$-regular sheaf is ample -- see 
\S7 for a nice application.  
\end{remark}

Let $Y$ be a smooth projective \emph{complex} variety whose Albanese map $a: Y \rightarrow 
{\rm Alb}(Y)$ is generically finite onto its image, in other words a variety of \emph {maximal Albanese dimension}. The well-known Generic Vanishing Theorem of Green-Lazarsfeld, and a related result of 
Ein-Lazarsfeld, imply in the present language the following: 

\begin{theorem}[\cite{gl} Theorem 1, \cite{el} (proof of) Theorem 3]\label{generic_vanishing}
If $Y$ is a variety of maximal Albanese dimension, then either $a_* \omega_Y$ 
is $M$-regular, or the Albanese image $a(Y)$ is ruled by subtori of ${\rm Alb}(Y)$. 
\end{theorem}

It is worth noting the conjecture of Koll\' ar (cf. \cite{kollar}, 17.9) saying that if $Y$ is of maximal Albanese dimension and of general type, then $\chi(\omega_Y) > 0$. By \cite{el} this is not always true, but the failure is accounted for by well-described special behavior:

\begin{corollary}[\cite{el}, Theorem 3]
If $Y$ is a smooth projective variety of maximal Albanese dimension with $\chi(\omega_Y)=0$, then the Albanese image $a(Y)$ is ruled by tori\footnote{The paper 
\cite{el} provides examples of such $Y$ of general type and with 
$\chi (\omega_Y) = 0$.} .
\end{corollary}
\begin{proof}
By the Theorem above, the only thing that needs to be said is that if 
$a_* \omega_Y$ is $M$-regular, then $\chi (\omega_Y) > 0$. But this follows 
immediately from Lemma \ref{nonvanishing} and Grauert-Riemenschneider vanishing. 
\end{proof}

Combining Theorem \ref{generic_vanishing} with the main $M$-regularity result,  Theorem \ref{mreg}, we obtain the following (cf. \cite{pp3} Theorem 5.1): 

\begin{theorem}\label{3can}
Let $Y$ be a smooth projective complex minimal variety of general type, of maximal Albanese dimension. If $a(Y)$ is not ruled by tori, then $|2K_Y|$ is basepoint-free and $|3K_Y|$ is very ample 
outside of the exceptional locus of $a$ (the union of the positive dimensional fibers). 
\end{theorem}

\begin{remark}
With extra work one can show that $|3K_Y|$ separates general points without the assumption 
that $Y$ is minimal. This was proved by Chen and Hacon in \cite{ch3}. 
They (cf. e.g. \cite{ch1}, \cite{ch2}) and Hacon-Pardini (cf. \cite{hp}) also proved many other very interesting results about the geometry of varieties of maximal Albanese dimension, based on the use  of generic vanishing theorems and the Fourier-Mukai transform. 
\end{remark}

In the case of adjoint line bundles, vanishing of higher cohomology is provided by the
Kawamata-Viehweg vanishing theorem, and is preserved via generically finite maps. 

\begin{lemma}\label{finite}
If $Y$ is a smooth projective variety of maximal Albanese dimension and $L$ is a big and nef line bundle on $Y$, then $a_* \OO(K_Y+ L)$ satisfies I.T. with index $0$ (so it is trivially $M$-regular). 
\end{lemma}
\begin{proof} 
By Kawamata-Viehweg vanishing 
$$H^i (\OO(K_Y + L)\otimes a^* \alpha) = 0, {\rm ~for~ all~} i> 0 {\rm ~and~ all~} 
\alpha\in {\rm Pic}^0 (Y).$$
On the other hand we have $R^i a_* \OO(K_Y + L) = 0$ for all $i >0$. This follows exactly as in the proof of Grauert-Riemenschneider vanishing (cf. \cite{positivity} \S4.3.B), which is the same statement 
for $K_Y$. By the degeneration of the Leray spectral sequence we have then
$$H^i (a_* \OO(K_Y+ L) \otimes \alpha) = 0, {\rm ~for~ all~} i> 0 {\rm ~and~ all~} 
\alpha\in {\rm Pic}^0 (Y).$$
\end{proof} 

Combining this in the usual way with Theorem \ref{mreg}, one obtains the following 
"generalized Lefschetz theorem" for adjoint line bundles. 

\begin{theorem}[cf. \cite{pp1}, Theorem 5.1] 
If $Y$ is a smooth projective variety of maximal Albanese dimension and $L$ is a big and nef line bundle on $Y$, then 
\begin{enumerate}
\item $|K_Y + L| \neq \emptyset$.\footnote{It is conjectured by Kawamata \cite{kawamata} that on every smooth projective variety $Y$, an adjoint line bundle $\OO(K_Y + L)$ with $L$ nef and big should have sections, as long as it is nef -- in our case this last property is automatic by Lemma \ref{nonvanishing}.}
\item $|2(K_Y + L)|$ is basepoint-free outside of the exceptional locus of $a$. 
\item $|3(K_Y + L)|$ is very ample outside of the exceptional locus of $a$. 
\end{enumerate}
\end{theorem}
\begin{proof} 
The first statement follows combining Lemma \ref{finite} and Lemma \ref{nonvanishing}. 
The second follows since by Proposition \ref{cont} the sheaf $a_* \OO(K_Y+ L)$ is continuously globally generated and, like global generation, this property holds if and only if it holds after push-forward by a finite map.
The third statement follows in a similar way, by proving the global generation of $\OO_Y(2(K_Y + L))\otimes \I_{\{y \}}$ for every $y \in Y$ outside of the exceptional locus.  
\end{proof}

This result was stated in \cite{pp1} only for finite Albanese maps -- however, as explained above, 
since Grauert-Reimenschneider vanishing also holds for adjoint-type bundles, the extension to generically finite Albanese maps is immediate. 

For other applications of $M$-regularity to questions related to linear series on 
special irregular varieties, cf. \cite{pp1} \S5 and \cite{pp3} \S5.

\section{Equations of curves and special varieties in Jacobians} 

Let $C$ be a smooth curve of genus 
$g\geq 3$, and denote by $J(C)$ the Jacobian of $C$. Let $\Theta$ be a theta
divisor on $J(C)$, and $C_d$ the $d$-th symmetric product of $C$, for $1\leq d\leq g-1$. Consider 
$$u_d :C_d\longrightarrow J(C)$$ 
to be an Abel-Jacobi mapping of the symmetric product (depending on the choice 
of a line bundle of degree $d$ on $C$), and denote by $W_d$ 
the image of $u_d$ in $J(C)$. The ideas involved in the study of regularity provide basic results on the 
cohomology of the ideal sheaves of these special subvarieties. 

\begin{theorem}[\cite{pp1} Theorem 4.1, \cite{pp3} Theorem 4.3]\label{w_d}
For any $1\leq d\leq g-1$, the twisted ideal sheaf 
$\I_{W_d}(2\Theta)$ satisfies I.T. with index $0$ 
(i.e. $\I_{W_d}$ is (strongly) $3$-$\Theta$-regular). 
\end{theorem}

\begin{corollary}\label{cubic_theta}
For any $1\leq d\leq g-1$, the ideal $\I_{W_d}$ is cut out by divisors
algebraically equivalent to $2\Theta$. 
Moreover, $\I_{W_d}$ is also cut out by divisors linearly equivalent to $3\Theta$.
\end{corollary}

The Corollary follows immediately from Theorem \ref{w_d}, combined with Theorem \ref{mreg} and Proposition \ref{cont} (cf. also Example \ref{basic} (3)).

Following the approach in \cite{pp1}, the proof of Theorem \ref{w_d} goes roughly as follows. 
From the standard cohomology sequence associated to the exact sequence
$$0\longrightarrow \I_{W_d}(2\Theta)\longrightarrow \OO_{J(C)}(2\Theta)
\longrightarrow \OO_{W_d}(2\Theta)\longrightarrow 0,$$
we see that the theorem follows as soon as we prove the following two statements, 
where this time for simplicity we denote by $\Theta$ any of the translates of the principal polarization:
\begin{enumerate}
\item $H^i(\OO_{W_d}(2\Theta))=0, ~\forall~i>0$.
\item The restriction map $H^0(\OO_{J(C)}(2\Theta))
\rightarrow H^0(\OO_{W_d}(2\Theta))$ is surjective.
\end{enumerate}

On the other hand, it is possible to see that the restriction of the principal polarization to $W_d$ satisfies the following properties: 

\begin{itemize} 
\item The sheaf $\OO_{W_d}(\Theta)$ is $M$-regular on $J(C)$.
\item We have $h^0(W_d, \OO_{W_d}(\Theta)\otimes \alpha)=1$ for 
$\alpha\in \widehat{J(C)}$ general.
\end{itemize} 

This is done by explicit computation (cf. \cite{pp1}, Proposition 4.4) -- for example one can 
see that the cohomological support loci are $V^i(\OO_{W_d}(\Theta)) \cong W_{g-d-1}$, for all 
$1\le d \le g-1$. The $M$-regularity of $\OO_{W_d}(\Theta)$ implies (1) above via Proposition
\ref{freg_vanishing}. 

For the second statement one uses Theorem \ref{amultiplication}. 
Note that for any open subset $U\in \widehat{J(C)}$ we have a commutative 
diagram as follows, where the vertical maps are the natural restrictions. 
$$\xymatrix{
\bigoplus_{\alpha \in U} H^0(\OO_{J(C)}(\Theta)\otimes \alpha)\otimes 
H^0(\OO_{J(C)}(\Theta)\otimes \alpha^{-1}) \ar[r] \ar[d] &
H^0(\OO_{J(C)}(2\Theta)) \ar[d] \\
\bigoplus_{\alpha \in U} H^0(\OO_{W_d}(\Theta)\otimes \alpha)\otimes 
H^0(\OO_{J(C)}(\Theta)\otimes \alpha^{-1}) \ar[r] & H^0(\OO_{W_d}(2\Theta)) } $$
Now Theorem \ref{amultiplication} says that the bottom horizontal map 
is surjective. On the other hand, by the second item above, we can choose 
the open set $U$ such that the left vertical map is an isomorphism \footnote{This 
means considering the translates $\Theta_\xi$ satisfying  
$h^0(\OO_{W_d}(\Theta_\xi))=1$, plus the fact that $W_d$ is not contained in $\Theta_\xi$.}. 
This in turn implies that the right vertical map is surjective, i.e. (2) above.

These vanishing results on Jacobians can also be approached via a study of Picard bundles 
associated to line bundles of high degree on $C$ -- the main tool is still the Fourier-Mukai 
transform (cf. \cite{pp3} \S 4). Finally Theorem \ref{w_d}, together with a Grothendick-Riemann-Roch computation, implies also a precise formula for  
the cohomology of pull-backs of $2\Theta$-line bundles to symmetric products, as conjectured 
by Oxbury and Pauly in \cite{oxbury}:

\begin{corollary}\label{op_conj}
In the notation above the following are true:
\newline
\noindent
(1) $h^0(C^d, u_d^*\OO(2\Theta))= \Sigma_{i=0}^d {g \choose i}$
\newline
\noindent
(2) $u_d^*:H^0(J(C),\OO(2\Theta))\rightarrow H^0(C^d, u_d^*\OO(2\Theta))$ is surjective.
\end{corollary}

\noindent
\textbf{A geometric Schottky-type criterion following Castelnuovo's Lemma.} One of the consequences of 
Theorem \ref{w_d} above is the following statement about points on Abel-Jacobi curves. 

\begin{corollary}
Let $\Gamma$ be an arbitrary set of 
$n\ge g+1$ points lying on a curve $C$ Abel-Jacobi embedded in its Jacobian. 
Then $\Gamma$ imposes precisely $g+1$ conditions on the linear series $|(2\Theta)_a|$ for $a\in J(C)$ general.
\end{corollary}

Indeed the restriction map from $2\theta$-functions on the Jacobian to the points of $\Gamma$ factors through the restriction to the curve: 
$$H^0 \OO_{J(C)} ((2\Theta)_a) \rightarrow H^0 \OO_C ((2\Theta)_a) \rightarrow 
H^0 \OO_{\Gamma} ((2\Theta)_a).$$
By Theorem \ref{w_d}, the first map is surjective, while the second one is clearly injective for $a$ general. 

On the other hand -- as a consequence of the Jacobi inversion theorem -- a 
collection of  $n\ge g+1$ general points on an Abel-Jacobi curve is in theta general 
position, an abelian varieties analogue of linear general position in $\PP^n$. 

\begin{definition}
A collection $\Gamma$ of $n\ge g+1$  distinct
points on a principally polarized abelian variety $(X, \Theta)$ is in 
\emph{theta general position} if  
for any $Y\subset \Gamma$ with $|Y|=g$ and for any $p\in
\Gamma-Y$ there is a theta-translate $\Theta_{\gamma}$ such
that $Y\subset \Theta_{\gamma}$ and $p\not\in
\Theta_{\gamma}$.
\end{definition}

It is easily checked that $g+1$ is the minimal number of conditions possibly 
imposed on the general linear series $|(2\Theta)_a|$ by points in theta general position. 
Thus general points on Abel-Jacobi curves impose the smallest possible number of 
conditions on ``$2\theta$-functions". The converse is also true, in the sense that this 
condition identifies precisely curves in their Jacobians. 

\begin{theorem}[\cite{pp4}, Theorem 5.2]\label{cs}
Let $(X,\Theta)$ be an irreducible principally polarized abelian
variety of dimension $g$, and let $\Gamma$ be a set of $n\geq g+2$
points on $X$ in theta general position, imposing only $g+1$
conditions on the linear series $|(2\Theta)_a|$ for $a\in X$
general. Then $(X, \Theta)$ is a canonically polarized Jacobian of
a curve $C$ and $\Gamma \subset C$ for a unique
Abel-Jacobi embedding $C\subset J(C)$.
\end{theorem}

This converse however does not use $M$-regularity, but rather results of Gunning and Welters on the 
existence of families of trisecants to the Kummer embedding. 
A related result has been obtained by Grushevsky in \cite{grushevsky} using the analytic 
theory of theta functions. 

Theorem \ref{cs} is an  abelian varieties analogue of the Castelnuovo Lemma (cf. \cite{gh} p.531) in projective space. Together with Theorem \ref{w_d}, it shows that the unique type of nondegenerate curve of minimal degree in projective space (the rational normal curve) behaves like the unique type of nondegenerate curve representing the minimal class  $\frac{\theta^{g-1}}{(g-1)!}$ in a ppav (the Abel-Jacobi curve). On the other hand, it is well known that varieties of minimal degree in $\PP^n$ are 
the only $2$-regular subvarieties in the sense of Castelnuovo-Mumford. Thus the following:

\begin{conjecture}[cf. \cite{pp4} \S2]
\emph{A non-degenerate subvariety $Y$ of an irreducible ppav $(X, \Theta )$
represents a minimal class $\frac{\theta^m}{m!}$, for any $m\ge 1$, if and only if $\I_Y(2\Theta)$ satisfies I.T. with index $0$ (i.e. $\I_Y$ is {\rm strongly} $3$-$\Theta$-regular).} 
\end{conjecture}

Subvarieties of ppav's representing the minimal classes $\frac{\theta^m}{m!}$ are believed (cf. \cite{debarre1}) to be -- up to translation and multiplication by $-1$ -- precisely the $W_d$'s in Jacobians, together with the Fano surface of lines in the intermediate Jacobian of the smooth cubic threefold \footnote{The vanishing in the conjecture is Theorem \ref{w_d} above for $W_d$'s, and is checked for the Fano surface in \cite{horing}.}. For more on the potential analogy with the situation in $\PP^n$ cf. \cite{pp4} \S2.

\section{Other applications}

\noindent
{\bf Maps to simple abelian varieties \cite{debarre2}.}
In \cite{debarre2} Debarre proves a general positivity result,  which improves substantially the observation in Lemma \ref{nonvanishing} (1)\footnote{The result is proved in \cite{debarre2} more generally for any coherent sheaf, with an interpretation of ampleness similar to that of 
nefness in Lemma \ref{nonvanishing}.}. 

\begin{proposition}\label{ampleness}
Every $M$-regular vector bundle on an abelian variety is ample. 
\end{proposition}

Let now $f: X \rightarrow Y$ be a finite surjective morphism of smooth projective complex varieties. The trace map ${\rm Tr}: f_* \OO_Y \rightarrow \OO_X$ splits, and this defines a complement vector bundle by 
$$ f_* \OO_Y  = \OO_X \oplus E_f^{\vee}.$$
Using the strong information provided by Generic Vanishing, Debarre shows in \emph{loc. cit.} the following: 

\begin{proposition}
If $Y$ is a simple abelian variety and the morphism  $f$ does not factor through a non-trivial isogeny, then the vector bundle $E_f$ is $M$-regular, hence ample. 
\end{proposition}

This is the analogue, for abelian varieties, of Lazarsfeld's result on finite covers of projective 
space, stating that if $Y = \PP^n$, then $E_f$ is ample -- it is also well-known by Lazarsfeld's work that the ampleness of $E_f$ has strong topological implications (cf. \cite{positivity}, Theorem 6.3.55 and Example 6.3.56). As  a consequence of the results above, Debarre obtains the nice application: 

\begin{theorem}[\cite{debarre2}, Theorem 1.1]
Let $f: X \rightarrow Y$ be a finite morphism of degree $d$ from a smooth connected projective variety to a simple abelian variety of dimension $n$. If $f$ does not factor through a nontrivial isogeny, then the induced morphism 
$$H^i(f, \CC): H^i (Y, \CC) \longrightarrow H^i(X, \CC)$$
is bijective for $i\le n-d+1$. 
\end{theorem}

\begin{remark}
This application illustrates once more the fact that $M$-regularity takes up the role played by Castelnuovo-Mumford regularity in projective space: Lazarsfeld's proof of the original result on finite maps to $\PP^n$ is based precisely on a study of the Castelnuovo-Mumford regularity of the bundle $E_f$. 
\end{remark}

\noindent
{\bf Verlinde bundles on moduli spaces of vector bundles on curves.}
Let $C$ be a smooth projective complex curve of genus $g\geq 2$. 
Let also $U_{C}(r,0)$ be the moduli space of semistable vector bundles of rank $r$ 
and degree $0$ on $C$ and $SU_{C}(r)$ the moduli space of semistable rank $r$ vector 
bundles with trivial determinant.
The \emph{Verlinde bundles} are push-forwards of  
pluritheta line bundles on $U_{C}(r,0)$ to the Jacobian of $C$ via the determinant 
map ${\rm det}:U_{C}(r,0)\rightarrow J(C)$:
$$E_{r,k}:={\rm det}_{*}\OO(k\Theta_{N}),$$
where $\Theta_{N}$ is the generalized theta divisor associated to a line bundle 
$N\in {\rm Pic}^{g-1}(C)$ (cf. \cite{popa} \S1).  Their fibers over $L\in J(C)$ are precisely the Verlinde 
vector spaces of level $k$ theta functions $H^{0}(SU_{C}(r,L),\mathcal{L}^{k})$, where $\mathcal{L}$ is 
the determinant line bundle on $SU_{C}(r)$.
The main technical result on these bundles is the following: 

\begin{proposition}[\cite{popa} Proposition 5.2 and Theorem 5.9]
The vector bundle $E_{r,k}$ is: 
\begin{itemize}
\item globally generated if and only if $k\ge r+1$.
\item normally generated if and only if $k\ge 2r+1$.
\end{itemize}
\end{proposition} 

The result is proved as a simple application of Theorem \ref{mreg}: the pull-back of 
$E_{r,k}$ by the map given by multiplication by $r$ on $J(C)$ is isomorphic to 
$\bigoplus \OO_{J(C)} (kr\Theta )$, and so it can be controlled cohomologically. 
It can also be seen as an example of the general behavior of semihomogeneous 
vector bundles on abelian varieties (cf. \cite{pp3}, \S6).
It is applied in \cite{popa}, and in a more refined form in \cite{ep}, in order to obtain 
effective results on generalized theta linear series on $U_C(r,0)$.

\providecommand{\bysame}{\leavevmode\hbox to3em{\hrulefill}\thinspace}

\end{document}